\documentclass[12pt]{article}
\usepackage{amsmath}
\usepackage{amssymb}
\usepackage{amsthm}

\newtheorem{thm}{Theorem}[section]

\newtheorem{lem}{Lemma}[section]
\newtheorem{defn}{Definition}[section]
\newtheorem{prop}{Proposition}[section]

\allowdisplaybreaks

\numberwithin{equation}{section}

\def\II{{\mathfrak J}}
\def\RR{{\mathcal R}}
\def\bbR{{\mathbb R}}
\def\bbZ{{\mathbb Z}}
\def\bbN{{\mathbb N}}

\def\bbE{{\bf E}}
\def\bbP{{\bf P}}
\def\TT{{\mathbb T_d}}
\def\Tk+{{\mathbb T_d^+}}
\def\KK{{\mathcal K}}
\def\SS{{\mathcal S_k}}
\def\SSU{{\mathcal S_{k+1}}}
\def\AA{{\mathcal A}}
\def\DD{{\mathcal D}}

\def\GG{{\mathcal G}}
\def\VV{{\cal V}}
\def\EE{{\cal E}}
\def\EEE{{\mathcal E}}
\def\XX{{\cal X}}

\def\Pp{{\mathfrak P}}

\def\eps{\varepsilon}

\let\phi=\varphi

\def\qed{\hfill $\square$}
\def\I{{\bf 1}}
\def\HHH{{\cal H}}

\def\dist{\mathop{\rm dist}}

\newcommand{\FM}[3]{\ensuremath{\mbox{\rm FM}(#1,#2,#3)}}
\newcommand{\gFM}[2]{\ensuremath{\mbox{\rm FM}(#1,#2)}}

\begin{document}

\title{Phase transition for the frog model\thanks{The authors
are thankful to
CAPES/PICD, CNPq (300226/97--7 and 300676/00--0)
 and FAPESP (97/12826--6) for
financial support.}}

\author{O.S.M.~Alves$^{~1}$, F.P.~Machado$^{~2}$,
S.Yu.~Popov$^{~2}$}

\maketitle {\footnotesize

\noindent $^{~1}$Instituto de Matem\'atica e Estat\'\i stica,
Universidade Federal de Goias, Campus Samambaia, Caixa Postal 131,
CEP 74001--970, Goi\^ania GO, Brasil

\noindent e-mail: \texttt{oswaldo@mat.ufg.br}

\smallskip

\noindent $^{~2}$Instituto de Matem\'atica e Estat\'\i stica,
Universidade de S\~ao Paulo, Rua do Mat\~ao 1010, CEP 05508--900,
S\~ao Paulo SP, Brasil

\noindent e-mails: \texttt{fmachado@ime.usp.br, popov@ime.usp.br} }

\begin{abstract}
We study a system of
simple random walks on graphs, known as {\it frog model}.
This model can be described as follows: There are active and
sleeping particles living on some graph~$\GG$. Each active particle
performs a simple random walk with discrete time and at each moment
it may disappear with probability~$1-p$. When an active particle
hits a sleeping particle, the latter becomes active.
Phase transition results and asymptotic
values for critical parameters are presented for $\bbZ^d$ and
regular trees.
\\[.3cm]
\end{abstract}

\section{Introduction and results}
\label{intro}
The subject of this paper is the so-called frog model with death,
which can be described as follows.
Initially there is a random number of particles
at each site of a graph $ \GG $.
A site of $ \GG $ is singled out and called its root.
All particles are sleeping at time zero, except for those
that might be placed at the root, which are active.
At each instant of time, each active particle may die with
probability $ (1-p) $. Once an active particle survives, it
jumps on some of its nearest neighbors, chosen with uniform probability,
performing a discrete time simple random walk (SRW) on $ \GG $.
Up to the time it dies, it activates all sleeping particles it hits
along its way. From the moment they are activated on, every such
particle starts to walk, performing exactly the same dynamics,
independent of everything else.

This model with $p=1$ (i.e., no death) is a discrete-time version
of the model for information spreading proposed by R.~Durrett
(1996, private communication), who also suggested the term
``frog model''.
The first published result on
this model is due to Telcs, Wormald~\cite{TW}, where it was
referred to as the ``egg model''. They proved that, starting
from the one-particle-per-site initial configuration,
almost surely the origin
will be visited infinitely often. Popov~\cite{Serguei}
proved that the same is true in dimension $d\geq 3$
for the initial configuration constructed as follows:
A sleeping particle (or ``egg'') is added into
each $x\neq 0$ with probability
$\alpha/\|x\|^2$, where~$\alpha$ is a large positive constant.
In Alves et al.~\cite{sofa} for the frog model with no death
it was proved that, starting from the one-particle-per-site
 initial configuration, the set of the
original positions of all active particles, rescaled by the
elapsed time, converges to a nonempty compact convex set.
In Alves et al.~\cite{Ravi} a similar result was obtained
for the case of random initial configuration.

The authors have learned about this version (i.e., with death) of
the frog model from I.~Benjamini.
 The goal of the present work is to study the asymptotic
dynamics of this particle system model,
regarding to the parameter~$ p $,
the graph where the random walks take place and the
initial distribution of particles.

Let us define the model in a formal way. We denote by $ \GG
=(\VV,\EE) $ an infinite connected
non-oriented graph of locally bounded
degree. Here $ \VV := \VV(\GG) $ is the set of vertices (sites) of
$ \GG $, and $ \EE := \EE(\GG) $ is the set of edges of $ \GG $.
Sites are said to be neighbors if they belong to a common edge.
The {\it degree\/} of a site $ x $ is the number
 of edges which have $ x $ as
an endpoint. A graph
is {\it locally bounded\/} if all its sites have finite degree.
Besides, a graph has {\it bounded degree\/}
if its maximum degree is finite.
The distance $ \dist(x,y) $ between sites~$ x $ and~$ y $ is the
minimal amount of edges that one must
pass in order to go from $ x $ to $ y $. Fix a site $ {\bf 0} \in \VV $
and call it the root of $ \GG $. With the usual abuse of notation,
by $ \bbZ^d $ we mean the graph with the vertex
set $ \bbZ^d $ and edge set $ \{
\langle(x^{(1)}; \ldots ; x^{(d)}),
(y^{(1)}; \ldots ; y^{(d)})\rangle :|x^{(1)}
- y^{(1)}| + \cdots + |x^{(1)} - y^{(d)}| = 1 \} $. Also, $ \TT,
\  d \ge 3 $, denotes the degree~$ d $ homogeneous tree.

Let $ \eta $ be a random variable taking values in $\bbN=\{0,1,2,
\dots\}$ such that $ \bbP [ \eta
\geq 1 ] > 0 $, and define
$ \gamma_j = \bbP[ \eta = j ] $. Let $ \{ \eta(x)
; x \in \VV \}, \{ (S^x_n(i))_{n \in \bbN}; i \in \{1,2,3, \dots
\}, x \in \VV \} $ and $ \{(\Xi_p^x(i)); i \in  \{1,2,3,
\dots \}, x \in \VV \} $ be independent sets of i.i.d.\
random variables defined as follows. For each $ x \in \VV$,
$ \eta(x) $ has the same law as $ \eta $, and gives the initial
number of particles at site~$ x $. If $ \eta(x)
\geq 1 $, then for each $ 0 < i \le \eta(x) $, $ (S^x_n(i))_{n \in
\bbN} $ is a discrete time SRW on $ \GG $ starting from $ x $
(it describes the trajectory of $i$-th particle from~$x$), and
$ \Xi_p^x(i) $, which stands for the lifetime of $i$-th particle from~$x$,
is a random variable whose law is given by $ \bbP[
\Xi_p^x(i) = k ] = (1- p) p^{k-1}$, $k=1,2,\ldots$,
where $p\in [0,1]$ is a fixed parameter.
%These sequences give the movement
%of each one of the $ \eta(x) $ particles
%placed originally at site~$ x $,
%starting to move at the instant they are activated (if
%that happens at all).
Thus, the $ i $-th particle at site $ x $ follows
the SRW  $ (S^x_n(i))_{n \in \bbN} $ and dies (disappears) $
\Xi_p^x(i) $ units of time after being activated.
For $x\neq y$ let
\[
 t(x,y) = \min_{1\leq i \leq \eta(x)} \min\{n<\Xi_p^x(i) : S_n^x(i)=y\}
\]
(clearly, $t(x,y)=\infty$ with positive probability). The moment when
all the particles in~$x$ are awakened is  defined as
\[
T(x) = \inf\{t(x_0,x_1)+\cdots+t(x_{m-1},x_m)\},
\]
where the infimum is over all finite sequences
$0=x_0, x_1, \ldots, x_{m-1},x_m=x$. Clearly, $T(x)=\infty$
means that the site~$x$ is never visited by active particles.

It is important to note that at the moment the
particle disappears, it is not able to activate other particles
(as first we decide whether the particle survives, and only after
that the particle that survived is allowed to jump).
Notice that there is no interaction between active particles, which
means that each active particle moves independently of everything
else. We denote by
\FM{\GG}{p}{\eta} the frog model on the graph $\GG$ with survival
parameter~$p$ and initial configuration ruled by~$\eta$.

Let us consider the following definition.

\begin{defn}
\label{d surv} A particular realization of the frog model\/ {\em
survives} if for any instant of time there is at least one active
particle. Otherwise, we say that it\/ {\em dies out}.
\end{defn}

Now we observe that $ \bbP [\FM{\GG}{p}{\eta}
\hbox{ survives}] $ is nondecreasing in $ p $ and define
\[
 p_c(\GG,\eta) := \inf \{ p :
\bbP [\FM{\GG}{p}{\eta} \hbox{ survives}] > 0 \}.
\]
 As usual, we say that \FM{\GG}{p}{\eta} exhibits
{\it phase transition\/} if
\[
 0 < p_c(\GG, \eta) < 1.
\]

Now we present two lower bounds on $p_c(\GG, \eta)$ which
can be obtained by a direct comparison with a Galton-Watson
branching process.
The next proposition shows that,
provided that $\bbE\eta<\infty$,
 for small enough~$p$ (depending on
$ \eta $) the frog model dies out almost surely on any graph.

\begin{prop}
\label{G}
If $\bbE\eta<\infty$, then for any graph $\GG$,
$p_c(\GG,\eta)\geq {(\bbE \eta +1)}^{-1}$.
\end{prop}

\noindent
{\it Proof.} Take $p\leq {(\bbE \eta +1)}^{-1}$.
The set of active particles in the frog
model is dominated by the population of the
following Galton-Watson branching
process. Each individual has a number of offspring distributed as
$ ( \eta +1) \xi $, where the random variable~$\xi$
is independent of~$\eta$,
and $\bbP[\xi = 1]= p = 1-\bbP[\xi=0]$. Therefore,
since the mean number of offspring by individual is $ (1+ \bbE
\eta )p$, the result follows by comparison with the Galton-Watson
branching process. \qed

\medskip

Next, again by comparison with Galton-Watson branching process, we give
another lower bound to $ p_c(\GG, \eta)$. This bound is better
than the one presented in Proposition~\ref{G} for bounded degree graphs.

\begin{prop}
\label{MD}
Suppose that $ \GG $ a graph of maximum degree~$ k $,
and $\bbE\eta<\infty$. Then
it is true that
\[
 p_c(\GG,\eta) \ge \frac{k}{1+(k-1)(\bbE \eta +1)}.
\]
\end{prop}

\noindent {\it Proof.} Consider a Galton-Watson branching process where
particles produce no offspring with probability $ 1 - p$, one
offspring with probability $ p/k $ and the random number
 $  \eta+1 $ of offspring with probability $ p(k-1)/k $.
  Observing that every site with
at least one active particle at time $ n > 0 $, has at least one
neighbor site whose original particle(s) has been activated prior to
time $n$, one gets that the frog model
is dominated by the Galton-Watson process
just defined. An elementary calculation shows that if $ p < k
(1+(k-1)(\bbE \eta +1))^{-1} $, the mean offspring
in the Galton-Watson process defined above
is less than one, therefore it dies out almost surely.
So, the same happens to the frog model. \qed

\medskip

Before going further, let us underline that in fact we
are dealing with percolation. Indeed,
 let
\[
 \RR_x^i = \{ S^x_n(i) : 0 \le n < \Xi_p^x(i) \} \subset \GG
\]
 be the ``virtual'' set of sites visited by the $i$-th particle
placed originally at~$ x $. The set $ \RR_x^i $ becomes ``real''  in
the case when~$ x $ is actually visited
(and thus all the sleeping particles from there
are activated). We define the (virtual)
range of site~$x$ by
\[
 \RR_x :=
\left\{
\begin{array}{ll}
\displaystyle \bigcup_{i=1}^{\eta(x)} \RR^i_x,
 & \mbox{ if } \eta(x) >0, \\
\vphantom{\displaystyle\sum^N}\{x\}, &  \mbox{ if } \eta(x) = 0.
\end{array}
\right.
\]
 Notice that the frog model survives if and only if there exists
an infinite sequence of distinct sites $ {\mathbf 0} = x_0, x_1,
x_2, \ldots $ such that, for all $j$,
\begin{equation}
\label{rangeperc}
x_{j+1} \in \RR_{x_j}.
\end{equation}
The last observation shows that
the extinction of the frog model
is equivalent to the finiteness of the cluster of~${\bf 0}$
in the following
oriented percolation model: from each site~$x$ the oriented edges
are drawn to all the sites of the set $\RR_x$.
 This approach is the key for the proof of most of the results
of this paper.

Next we state the main results of this paper. The
proofs are given in Section~\ref{proofs}.

\subsection{Extinction and survival of the process}

We begin by showing that, under mild conditions on the
initial number of particles, the process dies out a.s.\
(i.e., there is no percolation) in~$\bbZ$
for any $p<1$. From now on, $a\vee b$ stands for
$\max\{a,b\}$.

\begin{thm}
\label{Z1}
If $ \bbE \log (\eta \vee 1)< \infty $, then
 $ p_c(\bbZ,\eta) = 1 $.
\end{thm}

Next, we find sufficient conditions to guarantee that the
process becomes extinct for~$p$ small enough in $\bbZ^d$, $d\geq 2$,
and in $\TT$, $d\geq 3$.

\begin{thm}
\label{Tkext}
Suppose that there exists~$\delta>0$ such that $\bbE\eta^\delta<\infty$.
Then $p_c(\TT,\eta)>0$, i.e., the process on~$\TT$
dies out a.s.\ for~$p>0$ small enough.
\end{thm}

\begin{thm}
\label{Zdext}
Suppose that $\bbE (\log(\eta\vee 1))^d < \infty$. Then $p_c(\bbZ^d, \eta)>0$.
\end{thm}

Now, let us state the results related to the survival of the process.
 First, we show that for nontrivial~$\eta$
the frog model survives on $\bbZ^d$, $d\geq 2$, and on $\TT$, when
the parameter~$p$ is close enough to~$1$.

\begin{thm}
\label{Zd}
If $\bbP[ \eta \geq 1] > 0$, then $p_c(\bbZ^d,\eta)
< 1 $ for all $d \ge 2$.
\end{thm}

\begin{thm}
\label{Tk}
If $\bbP[ \eta \geq 1] > 0$, then $ p_c(\TT,\eta) < 1
$ for all $ d \ge 3 $.
\end{thm}

Now we state the counterpart of Theorem~\ref{Tkext}. Note
that Theorems~\ref{Tkext} and~\ref{Tksurv} give the
complete classification in~$\eta$ of the frog model on~$\TT$
from the point of view of positivity of $p_c(\TT,\eta)$.
\begin{thm}
\label{Tksurv}
If $\bbE\eta^\delta = \infty$ for any $\delta>0$, then
$p_c(\TT, \eta)=0$.
\end{thm}

Besides, we are able to show that, for any fixed~$d$, if~$\eta$
has a sufficiently heavy tail, then \FM{\bbZ^d}{p}{\eta}
survives with positive probability
for all values of $ p \in (0,1)$ (which would be the counterpart of
Theorem~\ref{Zdext}).
However, we do not state this result now, as in Section~\ref{sec_other}
we will prove a stronger result (cf.~Theorem~\ref{Zdrec}).

\subsection{Asymptotics for $p_c$}

 The following two theorems give asymptotic values for
critical parameters (compare with Propositions~\ref{G} and~\ref{MD})
for the case of~$\bbZ^d$ and regular trees.

\begin{thm}
\label{Tklim}
We have, for the case $\bbE\eta<\infty$,
\[
 \lim_{d \to \infty} p_c(\TT,\eta) = \frac{1}{1+\bbE
\eta}.
\]
\end{thm}

\begin{thm}
\label{Zdlim}
We have, for the case $\bbE\eta<\infty$,
\[
 \lim_{d \to \infty} p_c(\bbZ^d,\eta) = \frac{1}{1+\bbE \eta}.
\]
\end{thm}

\noindent
{\it Remark.} Observe that by truncating~$\eta$ and using a simple
coupling argument one gets that if $\bbE\eta=\infty$, then
\[
\lim_{d \to \infty} p_c(\TT,\eta) = \lim_{d \to \infty} p_c(\bbZ^d,\eta) =0.
\]

\medskip

Note that Theorems~\ref{Tklim} and~\ref{Zdlim} suggest
that there is some monotonicity of the critical probability
in dimension. Then, a natural question to ask is the following:
Is it true that $p_c(\bbZ^d,\eta) \geq p_c(\bbZ^{d+1},\eta)$
for all~$d$ (and can one substitute ``$\geq$'' by``$>$'')?
 In fact, there is a more general question: if
$\GG_1\subset\GG_2$, is it true
that $p_c(\GG_1,\eta)\geq p_c(\GG_2,\eta)$?
The last question has a trivial negative answer if we construct
$\GG_2$ from~$\GG_1$ by adding loops; if loops are not allowed,
then this question is open. Note that for percolation that
inequality is trivial; even the strict inequality can be proved
in a rather general situation, cf.\ Menshikov~\cite{Micha}.

\subsection{Other types of phase transition and generalizations}
\label{sec_other}
There are other types of phase transitions for this
model which may be of interest. For example, let~$p$ be such
 that $p_c(\GG,1)<p<1$ and~$\eta_q$ be a
0-1 random variable with $\bbP[\eta_q=1]=1-\bbP[\eta_q=0]=q$. Then,
the following result holds:

\begin{prop}
\label{pt_q}
There is a phase transition in~$q$, i.e.,
 \FM{\GG}{p}{\eta_q} dies out when~$q$ is small and survives
 when~$q$ is large.
\end{prop}

\medskip
\noindent
{\it Proof.}
First, note that \FM{\GG}{p}{\eta_q} is dominated by the following
Galton-Watson branching process: An individual has 0 offspring
with probability $1-p$, 1 offspring with probability $p(1-q)$, and
2 offspring with probability~$p q$. The mean offspring of this branching
process is $p(1+q)$, so \FM{\GG}{p}{\eta_q} dies out if $q<-1+1/p$.

Let us prove that \FM{\GG}{p}{\eta_q} survives when $p>p_c(\GG,1)$
and~$q$ is close enough to~1. Indeed, this model dominates a
model described in the following way: The process starts
  from the one-particle-per-site
initial configuration, and on each step active particles decide {\it twice\/}
whether to disappear, the first time  with probability~$1-q$, and
the second time  with probability~$1-p$. The latter model is in fact
\FM{\GG}{p q}{1}, so the model \FM{\GG}{p}{\eta_q} survives if
$q>p_c(\GG,1)/p$. \qed

\medskip

One may also be interested in the study of other types of critical
behaviour with respect to the parameter~$p$.
Consider the following
\begin{defn}
The model \FM{\GG}{p}{\eta} is called recurrent if
\[
\bbP[{\bf 0}\mbox{ is hit infinitely
often in }\FM{\GG}{p}{\eta}]>0.
\]
 Otherwise, the model is called transient.
\end{defn}
Note that,
even in the case when
a single SRW on~$\GG$ is transient,
it is still reasonable to expect that the frog model with~$p=1$
is recurrent.
 For example, for the model \FM{\bbZ^d}{1}{1}
the recurrence was established in~\cite{TW}.
However, establishing the recurrence
property in this case is nontrivial;
it is still unclear to us whether \FM{\TT}{1}{1} is recurrent.
Now, denote
\[
p_u(\GG,\eta) = \inf\{p:\bbP[{\bf 0}\mbox{ is hit infinitely
often in }\FM{\GG}{p}{\eta}]>0\}
\]
(here, by definition, $\inf \emptyset = 1$);
clearly, $p_u(\GG,\eta) \geq p_c(\GG,\eta)$ for any $\GG$ and $\eta$.
Now, we are interested in studying the existence of phase
transition with respect to~$p_u$.

First, we discuss some situations when the model is transient for
any~$p$, except possibly the case $p=1$.
\begin{thm}
\label{Tktrans}
Suppose that $\bbE\eta^\eps < \infty$ for any~$0<\eps<1$.
Then $p_u(\TT,\eta)=1$.
\end{thm}

\begin{thm}
\label{Zdtrans}
Suppose that $\bbE[\log(\eta\vee 1)]^d < \infty$.
Then $p_u(\bbZ^d,\eta)=1$.
\end{thm}

The next two theorems give sufficient conditions to have $p_u<1$
on trees and on $\bbZ^d$.

\begin{thm}
\label{Tkrec}
Suppose that there exists $\beta < \frac{\log (d-1)}{2\log d}$
such that
\begin{equation}
\label{eq_Tkrec}
\bbP[\eta \geq n] \geq \frac{1}{n^\beta}
\end{equation}
for all~$n$ large enough. Then $p_u(\TT, \eta)<1$.
\end{thm}

\begin{thm}
\label{Zdrec}
Suppose that there exist $\beta<d$ such that
\begin{equation}
\label{eq_Zdrec}
\bbP[\eta \geq n] \geq \frac{1}{(\log n)^\beta}
\end{equation}
for all~$n$ large enough. Then $p_u(\bbZ^d, \eta)=0$.
\end{thm}

\noindent
{\it Remark.}
It is possible to see that,
if $\bbP[\eta \geq n] \geq (\log n)^{-\beta}$
for some $\beta<1$ and all~$n$ large enough, then \FM{\GG}{p}{\eta}
is recurrent on {\it any\/} infinite
connected graph~$\GG$ of bounded degree.
 Indeed, for arbitrary
graph~$\GG$ of bounded degree we do the following: First, fix a
subgraph~$\GG_1$ of~$\GG$, which is isomorphic to $\bbZ_+$. If~$k_0$
is the maximal degree of~$\GG$, it is easy to see that \FM{\GG}{p}{\eta}
dominates \FM{\GG_1}{p/k_0}{\eta} (if a particle wants to leave~$\GG_1$,
we just erase this particle). Then, we just apply Theorem~\ref{Zdrec}
for the case of~$\GG=\bbZ$ (from the proof of Theorem~\ref{Zdrec}
one gets that the argument for the case of~$\bbZ$ also works
for~$\bbZ_+$).
\medskip

 Theorems~\ref{Tkrec} and~\ref{Zdrec}  give sufficient
conditions on the tail of the distribution of~$\eta$
for the process to be recurrent when~$p<1$.
On the other hand, Theorems~\ref{Tktrans} and~\ref{Zdtrans} show
that for the one-particle-per-site initial
configuration the process is not recurrent even the parameter~$p<1$
is very close to~$1$. Being the model with one-particle-per-site initial
configuration the most natural example one has
to hand, a natural question is raised: What can be done
(i.e., how can one modify the model) to make the model
recurrent without augmenting the initial configuration.
Notice that, by definition, in our model the lifetime of active particles
is geometrically distributed.
In order to find answers to that question, we are going to change this
and study the situation when the lifetime has another
distribution, possibly more heavy-tailed one.

Let $\Xi$ be any nonnegative integer-valued
random variable. From this moment on we study the frog model
on~$\GG$ with one-particle-per-site initial
configuration, and the lifetimes of particles after activation
are i.i.d.\ random variables $(\Xi^x, x\in\GG)$ having the same
law as~$\Xi$. This model will be called \gFM{\GG}{\Xi}.

\begin{thm}
\label{chi_trans}
Suppose that one of the following alternatives holds:
\begin{itemize}
\item $\GG=\bbZ$ and $\bbE\sqrt{\Xi}<\infty$,
\item $\GG=\bbZ^2$ and $\bbE\frac{\Xi}{\log(\Xi\vee 2)}<\infty$,
\item $\GG=\TT$ or $\bbZ^d, d\geq 3$ and $\bbE\,\Xi<\infty$.
\end{itemize}
Then \gFM{\GG}{\Xi} is transient.
\end{thm}

\begin{thm}
\label{chi_rec}
There exists a sequence of positive numbers
$\beta_1, \beta_2,\beta_3,\ldots$ such that if for all~$n$ large
enough one of the following alternatives holds
\begin{itemize}
\item $d=1$ and $\bbP[\Xi \geq n^2] \geq \beta_1 n^{-1}
 \log\log n$,
\item $d=2$ and $\bbP[\Xi \geq n^2] \geq \beta_2 n^{-2}
 (\log n)^2$,
\item $d\geq 3$ and $\bbP[\Xi \geq n^2] \geq \beta_d n^{-2}
 \log n$,
\end{itemize}
then \gFM{\bbZ^d}{\Xi} is recurrent.
\end{thm}

In fact, results of Popov~\cite{Serguei}
suggest that the following is true:
For $d\geq 3$ there exist $\tilde\alpha_0(d), \tilde\alpha_1(d)$ such that
if $\bbP[\Xi \geq n^2] \leq \tilde\alpha_0(d) n^{-2}$
for all~$n$ large enough,
then \gFM{\bbZ^d}{\Xi} is
transient, and if $\bbP[\Xi \geq n^2] \geq \tilde\alpha_1(d) n^{-2}$
for all~$n$ large enough,
then \gFM{\bbZ^d}{\Xi} is recurrent.
The heuristic explanation for this is
as follows. The particle originally in~$x$ has a good chance
(i.e., comparable with $\|x\|_2^{2-d}$,
where $\|\cdot\|_2$ is the Euclidean norm)
of ever getting to the origin
only if it lives at least of order $\|x\|_2^2$
units of time (cf.\ Lemma~\ref{pofvisit}
below), so one may expect
that \gFM{\bbZ^d}{\Xi} behaves roughly as the frog model with
infinite lifetime of the particles and the initial configuration
of sleeping particles constructed as follows: we add a sleeping particle
into~$x$ with probability $h(x):=\bbP[\Xi \geq n^2]$, and add nothing
with probability $1-h(x)$. For the case when $h(x)\simeq \alpha/\|x\|_2^2$
the latter model was studied in~\cite{Serguei} and it was proved that
it is recurrent when~$\alpha$ is large and transient when~$\alpha$
is small (note that the transience also can be proved by dominating the
frog model by a branching random walk,
cf.\ e.g.\ den Hollander et al.~\cite{HMP}).
However, turning this heuristics into a rigorous proof
is presently beyond our reach.

%H\"aggstr\"om~\cite{H} recently studied the connection of existence
%of phase transitions for different percolation and particle
%system models on the same graph. For the frog model, one may
%also ask a similar question. Denote by $p_{perc}(\GG)$ the critical
%probability for the site percolation on the graph~$\GG$.
%Now, is it true that
%\[
%p_{perc}(\GG)\in (0,1) \Longleftrightarrow p_c(\GG,1)\in (0,1)
%\]
%for any graph $\GG$ of bounded degree?
%Note that the backward implication does not hold for graphs
%of unbounded degree. As an example, consider $\GG={\mathbb T}_!$,
%where ${\mathbb T}_!$ is a tree constructed in such a way that
%for any $v\in{\mathbb T}_!$ the degree of~$v$ is equal to
%$\dist({\bf 0},v)+2$. Then it is not difficult to see that
%$p_{perc}({\mathbb T}_!)=0$, while $p_{c}({\mathbb T}_!,1)=1/2$.
%Note also that for quasi-transitive graphs on which the SRW is
%recurrent the forward implication can be proved quite
%analogously to the proof of Theorem~\ref{Zd} (the case of $d=2$).

\section{Proofs}
\label{proofs}

\subsection{Preliminaries}
\label{Prelim}
Here we state a few basic facts which will be necessary later
in the Sections~\ref{s_extsurv}, \ref{s_asympt}, and~\ref{s_rectrans}.

For $0\leq p\leq 1$ and integer numbers $k,i\geq 1$ denote
$\Phi(i,k,p)=1-(1-p^k)^i$ and
$\hat k(i,p)=\lfloor\log i/\log(1/p)\rfloor$,
where $\lfloor x\rfloor$ stands for the largest integer
which is less than or equal to~$x$. The following fact can
be easily obtained by using elementary calculus and is stated
without proof.
\begin{lem}
\label{analiz}
There exist constants $\hat \beta_1, \hat \beta_2, \hat \beta_3$
such that for all $i,p$
\[
\hat \beta_1 \leq \Phi(i,k,p) \leq 1
\]
for $k\leq \hat k(i,p)$ and
\[
\hat \beta_2 p^{k-\hat k(i,p)} \leq \Phi(i,k,p)
                     \leq \hat \beta_3 p^{k-\hat k(i,p)-1}
\]
for $k\geq \hat k(i,p)+1$.
\end{lem}

In the sequel we will make use of the following large
deviation result:
\begin{lem} [Shiryaev~\cite{Shir},  p.~68.]
\label{Shiryaev} Let $ \{X_i, i \ge 1\} $ be i.i.d.\ random variables with
$\bbP[X_i=1] = p$ and $\bbP[X_i=0] = 1-p$. Then for any
$0<p<a<1$ and for any $N \geq 1$ we have
\begin{equation}
\label{shi}
\bbP\Big[\frac{1}{N}\sum_{i=1}^N X_i \geq a\Big]
\leq \exp\{-N H(a,p)\},
\end{equation}
where
\[
 H(a,p) = a \log\frac{a}{p} + (1-a)\log\frac{1-a}{1-p} > 0.
\]
If $0<a<p<1$, then~(\ref{shi}) holds with $\bbP[N^{-1}\sum_{i=1}^N X_i\leq a]$
in the left-hand side.
\end{lem}

In order to prove Theorems~\ref{Zd} and \ref{Zdlim} we need some
auxiliary fact about projections of percolation models.
Let $\HHH^d:=\{H^d_x:x \in \bbZ^d\}$
be a collection of random sets such that $x \in H^d_x$ and the sets
$H^d_x - x$, $x \in \bbZ^d$, are i.i.d.
By using the sets of $\HHH^d$, define an oriented percolation
process on $\bbZ^d$ analogously to what was done for the frog model
(compare with (\ref{rangeperc})):
If there is an infinite sequence
of distinct sites, ${\bf 0}=x_0, x_1, x_2, \dots$,
such that $x_{j+1} \in H^d_{x_j}$ for
all $j=0,1,2, \dots$, we say that the cluster
of the origin is infinite or, equivalently, that $\HHH^d$ survives.

Let $\Lambda:=\{x \in \bbZ^d : x^{(k)}=0 \hbox{ for } k \ge 3 \}
\subset \bbZ^d$
be a copy of $\bbZ^2$ immersed into $\bbZ^d$ and
$\Pp : \bbZ^d \to \Lambda$ be the projection on the first two
coordinates. Let $\HHH^2:=\{H^2_x:x \in \Lambda\}$ be a collection of
random sets such that $x \in H^2_x$ and the sets
$H^2_x - x$, $ x \in \Lambda $, are i.i.d. Analogously,
one defines the percolation of the collection~$\HHH^2$.

\begin{lem}
\label{H}
Suppose that there is a coupling of
$H^d_{\bf 0}$ and $ H^2_{\bf 0}$ such that
\begin{equation}
\label{aga}
\Pp(H^d_{\bf 0}) \supset H^2_{\bf 0},
\end{equation}
i.e., the projection of $H^d_{\bf 0}$ dominates $ H^2_{\bf 0}$.
 Then
\[
\bbP[\HHH^d \mbox{ survives}] \geq \bbP[\HHH^2 \mbox{ survives}].
\]
\end{lem}

\medskip
\noindent
{\it Proof.}
The proof of this fact is standard and can be done by
carefully growing the cluster in~$\Lambda$ step by step,
and comparing it with the corresponding process in~$\bbZ^d$.
See e.g.\ Menshikov~\cite{M} for details.
\qed

\medskip

Let ${\hat q}( n, x )$ be the probability that a SRW (starting
from the origin) hits~$x$ until the moment~$n$.
The following fact about hitting probabilities of SRW
is proved in~\cite{sofa}, Theorem~2.2 (except for the case $d=1$).
\begin{lem}
\label{pofvisit}
\begin{itemize}
\item If $ d = 1, x \neq 0 $ and $  n \geq \|x\|_2^2 $,
then there exists a number $ w_1 > 0 $ such that
\begin{equation}
{\hat q}( n, x ) \ge w_1 .
\end{equation}
\item If $ d = 2, x \neq 0 $ and $  n \geq \|x\|_2^2 $,
then there exists a number $ w_2 > 0 $ such that
\begin{equation}
{\hat q}( n, x ) \ge \frac{w_2}{\log \|x\|_2} .
\end{equation}
\item Suppose that $ d \ge 3,  x \not= 0 $ and $ n \ge \|x\|_2^2 $.
Then there exists a collection of positive numbers
$ w_d > 0 $, $d\geq 3$, such that
\begin{equation}
{\hat q}( n, x) \ge \frac{w_d}{\| x \|_2^{d-2}} .
\end{equation}
\end{itemize}
\end{lem}

\noindent {\it Proof.}
To keep the paper self-contained, we give the proof
of this fact.
Let ${\hat p}_n(x)$ be the probability that
the SRW is in~$x$ at time~$n$, and $\tau_x$ be the moment
of the first hitting of~$x$. Also, denote
by~$G_n(x) = \sum_{k=0}^n {\hat p}_k(x)$
the mean number of visits to~$x$ until the moment~$n$ ($G_n(x)$ is
usually called {\it Green's function\/}).

Suppose without loss
of generality that $\|x\|_2^2 \leq n < \|x\|_2^2+1$. Observe that
\begin{eqnarray*}
G_n(x) &=& \sum_{j=0}^n {\hat p}_j(x) = \sum_{j=0}^n
\sum_{k=0}^j {\hat p}_k(0) \bbP[\tau_x=j-k] \\
&=& \sum_{k=0}^n {\hat p}_k(0) {\hat q}(n-k,x)
\le {\hat q}(n,x)G_n(0).
\end{eqnarray*}
So
\[
 {\hat q}(n,x) \ge \frac{G_n(x)}{G_n(0)} \ge
\left\{
\begin{array}{ll}
\displaystyle
\frac{\sum_{j={\lfloor n/2 \rfloor}}^n
{\hat p}_j(x)}{\sum_{j=0}^n{\hat p}_j(0)}, & d=1,2, \\
\vphantom{\sum^{\int^N}}
(G_\infty(0))^{-1} \sum_{j={\lfloor n/2 \rfloor}}^n
{\hat p}_j(x), & d \geq 3.
\end{array}
\right.
\]
 Using Theorem~1.2.1 of~\cite{Lawler}, after some elementary computations
we finish the proof. \qed

\subsection{Extinction and survival}
\label{s_extsurv}
{\it Proof of Theorem~\ref{Z1}.} Notice that, for any graph~$\GG$
 and all $ x \neq y \in \GG $, the following inequality holds:
\begin{equation}
\label{ur*}
 \bbP[ y \in \RR^1_x ] \leq p^{\dist(x,y)}.
\end{equation}
Clearly, for a fixed $y\in\bbZ$, we have
\[
\bbP[ y \notin \RR_x \hbox{ for all } x \neq y] =
   \prod_{x:x\neq y} (1-\bbP[y\in \RR_x]),
\]
so the left-hand side of the above display is positive if and
only if $\sum_{x:x\neq y}\bbP[y\in \RR_x]<\infty$. Now, by
using~(\ref{ur*}) and Lemma~\ref{analiz}, for some $C_1, C_2>0$
(depending only on~$p$) one gets
\begin{eqnarray*}
\sum_{x:x\neq y}\bbP[y\in \RR_x] &=& 2\sum_{x\geq 1}\bbP[{\bf 0}\in \RR_x] \\
&\leq & 2\sum_{k=1}^\infty\sum_{i=1}^\infty\gamma_i\Phi(i,k,p)\\
&=& 2\sum_{i=1}^\infty\gamma_i \Big(\sum_{k\leq \hat k(i,p)}\Phi(i,k,p)+
         \sum_{k\geq \hat k(i,p)+1}\Phi(i,k,p)  \Big) \\
&\leq & 2\sum_{i=1}^\infty \gamma_i (C_1\log i + C_2)<\infty.
\end{eqnarray*}
Thus, $\bbP[ y \notin \RR_x \hbox{ for all } x \neq y]>0$, so,
 by the ergodic theorem there is an infinite
sequence of sites $ \cdots < y_{-1} < y_0 < y_1 < \cdots $ such
that for all~$ i $, $y_i \notin \RR_x \hbox{ for all } x \neq y_i $.
 Therefore, for almost every
realization there is an infinite number of blocks of sites
without ``communication'' with its exterior, which prevents the active
particles to spread out. The result follows. \qed

\medskip
\noindent
{\it Proof of Theorems~\ref{Tkext} and~\ref{Zdext}.}
For $\GG=\bbZ^d$ or $\TT$ denote $s_k(\GG)=|\{y\in\GG : \dist(x,y)=k\}|$
(note that the right-hand side does not depend on the choice of the site~$x$).
Using~(\ref{ur*}) and Lemma~\ref{analiz},
one gets that for some positive constants
$C_1,C_2,C_3,C_4$
\begin{eqnarray*}
\bbE |\RR_x\setminus\{x\}| &=& \sum_{y:y\neq x} \bbP[y\in\RR_x]\\
 &=& \sum_{y:y\neq x}\sum_{i=1}^\infty \gamma_i
            \bbP[y\in\RR_x \mid \eta(x)=i]\\
 &\leq & \sum_{i=1}^\infty \gamma_i\sum_{k=1}^\infty s_k(\GG)
                    \Phi(i,k,p)\\
 &\leq & \sum_{i=1}^\infty \gamma_i \Big(\sum_{k=1}^{\hat k(i,p)}
     s_k(\GG) + \sum_{k=1}^\infty \hat \beta_3 s_k(\GG)p^{k-1} \Big)\\
 &\leq &   \left\{ \begin{array}{ll}
                     C_1 \displaystyle \sum_{i=1}^\infty \gamma_i
                     i^{\frac{\log(d-1)}{\log(1/p)}} + C_2, &
                     \GG=\TT,\\
                     C_3 \displaystyle \sum_{i=1}^\infty \gamma_i
                     \Big(\frac{\log i}{\log(1/p)}\Big)^d + C_4, &
                     \GG=\bbZ^d,
                   \end{array}
           \right.
\end{eqnarray*}
which is finite for all~$p<1$ in the case $\GG=\bbZ^d$ and
for~$p$ small enough in the case of $\GG=\TT$. Now, as for some~$p_0>0$
(which may depend on the graph~$\GG$) the convergence is uniform
in~$[0,p_0]$, there exists small enough~$p$
(which depends on~$\GG$) such that
$\bbE|\RR_x\setminus\{x\}|<1$ for \FM{\GG}{p}{\eta},
so one gets the proof by means
of domination by a subcritical branching process.
\qed

\medskip

In order to prove Theorems~\ref{Zd} and~\ref{Tk}
it is enough to show that, for~$ p $ large enough, the frog model
survives with positive probability in~$ \bbZ^d $, $ d \ge 2 $,
and in~$ \TT $, $d\geq 3$. Let us define the modified initial configuration
$ \eta' $ by
\[
  \eta'(x) = {\mathbf 1}_{\{\eta(x) \ge 1\}}.
\]
 Since \FM{\GG}{p}{\eta} dominates
\FM{\GG}{p}{\eta'}, without loss of
generality we prove Theorems~\ref{Zd} and~\ref{Tk} assuming that
the initial configuration is given by $ \{\eta'(x): x \in
\GG\} $.

\medskip
\noindent
{\it Proof of Theorem~\ref{Zd}.}
We start by considering the two dimensional frog model
\FM{\bbZ^2}{p}{\eta} which is equivalent to \FM{\Lambda}{p}{\eta},
since $\Lambda$ is a copy of $\bbZ^2$ (recall the notation~$\Lambda$ from
Section~\ref{Prelim}, it was introduced just before Lemma~\ref{H}).
It is a well-known fact that the
two-dimensional SRW (no death) is recurrent.
Then, given $ N \in \bbN $ and
assuming $ \eta' ({\mathbf 0}) = 1 $,
for sufficiently large $p=p(N)$, the probability that the first active
particle hits all the sites in the square $ [-2N, 2N]^2 \cap \bbZ^2 $
before dying can be made arbitrarily large. Besides, the
probability that there is a site~$ x \in [0, N)^2 \cap
\bbZ^2 $ such that $ \eta'(x) = 1 $, also can be made arbitrarily
large by means of increasing~$ N $.

Let us define now a two-dimensional percolation
process in the following way.
 Divide $ \bbZ^2 $ into disjoint squares of
side $ N $, i.e., write
\[
 \bbZ^2 = \bigcup_{(r; k)\in \bbZ^2} Q(r,k),
\]
 where $ Q(r,k) = (r N; k N) + [0, N)^2 \cap \bbZ^2 $.
Declare $ Q(r, k) $ open if the following happens (and closed
otherwise)
\[
 \bigcup_{x \in Q(r,k)} ( \{\eta'(x)=1\} \cap
\{\RR^1_x\supseteq([-2N,2N]^2+x)\cap\bbZ^2\})
\neq \emptyset.
\]
Observe that the events $\{Q(r,k) \mbox{ is open}\}$, $(r,k)\in\bbZ^2$,
are independent.
Notice that the frog model dominates this percolation process in
the sense that if there is percolation then the frog model
survives. It is not difficult to see that by suitably choosing
 $ N $ and $ p $ it is possible to make
 $\bbP[Q(r,k) \mbox{ is open}]$ arbitrarily close to~$1$, so the
percolation process can be made supercritical,
and thus the result follows for $ \bbZ^2 $.

\medskip

Now, by using Lemma~\ref{H}, we give the proof for dimensions
$d \ge 3$.  Let $\{\RR_y(p,d):y \in \bbZ^d\}$ be the collection of
the ranges for the $d$ dimensional frog model.
For the moment we write $\RR_y(p,d)$ instead of
$\RR_y$ to keep track of the dimension and the
survival parameter. Analogously, let
$\{\RR_y(p,2):y \in \Lambda\}$ be the collection of
the ranges for the two dimensional frog model
immersed in $\bbZ^d$. Notice that
$\Pp (\RR_x(p,d)) \subset \Lambda$ is
distributed as $\RR_{\Pp(x)}(p^{\prime}, 2)$ for
$p^{\prime}=2p/(d(1-p) + 2p)$, where, as before,
$\Pp: \bbZ^d \to \Lambda$ is the projection on
the first two coordinates. Since the fact $p^{\prime}<1$
implies that $p<1$, and, as we just have proven,
 $p_c(\bbZ^2,\eta)<1$ when $\bbE \eta < \infty$,
by using Lemma~{\ref{H}} we
finish the proof of Theorem~\ref{Zd}. \qed

\medskip
\noindent
{\it Proof of Theorem~\ref{Tk}.}
As in the previous theorem, we work with~$\eta'$ instead of~$\eta$.
In order to prove the result, we need some additional notation.
Notice that for any $a\in \TT$ there is a unique path connecting~$a$
to~${\bf 0}$; we write $a\geq b$ if~$b$ belongs to that path.
For $a\neq {\bf 0}$ denote $\Tk+(a)=\{b\in \TT : b\geq a\}$.
Fix an arbitrary site~$a_0$ adjacent to the root and let
$\Tk+=\TT\setminus\Tk+(a_0)$.
For any $A\subset\TT$ let us define the {\it external boundary\/}
$\partial_e(A)$ in the following way:
\[
\partial_e(A) = A\setminus \{a\in A : \mbox{ there exists
$b\in A$ such that } b>a\}.
\]
A useful fact is that if $A,B$ are finite and $A\subset B$,
then $|\partial_e(A)|\leq |\partial_e(B)|$.
Now, denote by $W_t$ the set of sites visited
until time~$t$ by a SRW (no death) in~$\TT$ starting from~${\bf 0}$.
Note that
\begin{itemize}
\item as SRW on tree is transient, one gets that with positive
probability $W_t\subset\Tk+$ for all~$t$;
\item $|\partial_e(W_t)|$ is a nondecreasing sequence, and, moreover,
it is not difficult to show that $|\partial_e(W_t)|\to\infty$ a.s.\
as $t\to\infty$.
\end{itemize}
The above facts show that for~$p$ large enough
\begin{equation}
\label{fgh}
\bbE |\partial_e(\RR_{{\bf 0}})\cap \Tk+| > \frac{1}{\bbP[\eta'>0]}.
\end{equation}
Now, all the initially sleeping particles
in $\partial_e(\RR_{{\bf 0}})\cap \Tk+$
are viewed as the offspring of the first particle.
By using~(\ref{fgh}) together with
the fact that for any $x,y\in\partial_e(\RR_{{\bf 0}})$,
$x\neq y$, we have $\Tk+(x) \cap \Tk+(y)=\emptyset$,
one gets that the frog model dominates a Galton-Watson branching
process with mean offspring greater than~$1$, thus concluding the
proof of Theorem~\ref{Tk}. \qed

\medskip
\noindent
{\it Proof of Theorem~\ref{Tksurv}.}
First, note the following fact: For any graph~$\GG$ with
maximal degree~$d$, it is true that
\begin{equation}
\label{ur**}
\bbP[y\in\RR^1_x] \geq \Big(\frac{p}{d}\Big)^{\dist(x,y)}.
\end{equation}
Keeping the notation $\Tk+(x)$ from the proof of Theorem~\ref{Tk},
denote
\[
L_k(x) = \{y\in\Tk+(x) : \dist(x,y)=k\}.
\]
Using~(\ref{ur**}) and Lemma~\ref{analiz}, we have
\begin{eqnarray*}
\bbE|\RR_x\cap L_{\hat k(i,p/d)}(x)| &\geq& \sum_{i=1}^\infty \gamma_i
  (d-1)^{\hat k(i,p/d)} \Phi(i,\hat k(i,p/d),p/d) \\
  &\geq& \hat\beta_1 \sum_{i=1}^\infty \gamma_i
                      i^{\frac{\log(d-1)}{\log(d/p)}} = \infty,
\end{eqnarray*}
so, by dominating a supercritical branching process by the frog
model, one concludes the proof.
\qed

\subsection{Asymptotics for the critical parameter}
\label{s_asympt}
{\it Proof of Theorem~\ref{Tklim}.} By Proposition~\ref{G},
$p_c(\TT,\eta) \ge {(1+\bbE\eta)}^{-1}$. So, it is enough to show that,
fixing $ p > {(1+\bbE\eta)}^{-1}$, the model \FM{\TT}{p}{\eta}
survives for~$ d $ large enough.

Let us define
\begin{equation}
\label{etas}
\eta^{(s)}=
\begin{cases}
\eta, &\hbox{if } \eta \le s,\\
 0, &\hbox{if } \eta > s.
\end{cases}
\end{equation}
By the monotone convergence theorem it follows that $ \bbE\eta^{(s)}
\to \bbE\eta $ as $ s \to \infty$, so, if $ p >
{(1+\bbE\eta)}^{-1} $, then it is possible to choose $ s $ large enough
so that $ p > {(1+\bbE\eta^{(s)})}^{-1}$. Fixed $ s $ and $ p$, notice
that \FM{\TT}{p}{\eta} dominates \FM{\TT}{p}{\eta^{(s)}}
 in the sense that if the latter survives
with positive probability, the same happens to the former.
Therefore, it is enough to show that if $ p >{(1+\bbE\eta^{(s)})}^{-1}$,
 then \FM{\TT}{p}{\eta^{(s)}} survives for~$ d $
sufficiently large.

Let $\xi_n$ be the set of active particles of \FM{\TT}{p}{\eta^{(s)}},
which are at level~$n$ (i.e., at distance~$n$ from the root) at time~$n$.
Next we to prove that
there exists a discrete time supercritical branching process,
which is dominated by $\xi_n$. We do
this by constructing an auxiliary process $\tilde\xi_n\subset\xi_n$.
 First of all, initially
the particle(s) in ${\bf 0}$ belong(s) to $\tilde\xi_0$. In general,
the process $\tilde\xi_n$ is constructed by the following rules.
If at time $n-1$ the set of particles $\tilde\xi_{n-1}$,
which lives on the level~$n-1$, is constructed,
then at time~$n$  the set of particles $\tilde\xi_n$
(which all are at level~$n$) is constructed in the following way.
Introduce some ordering of the particles of~$\tilde\xi_{n-1}$,
they will be allowed to jump according to that order.
Now, if the current particle survives, then
\begin{itemize}
\item if the particle jumps to some site of
level~$n$ and does not encounter
any particles that already belong to $\tilde\xi_n$
 there, then this particle as well as all the particles possibly
activated by it enter to $\tilde\xi_n$;
\item otherwise it is deleted.
\end{itemize}

The particles of $\tilde\xi_{n+1}$ activated by some particle from
$\tilde\xi_n$ are considered as the offspring of that particle;
note that, due to the asynchronous construction of the process
$\tilde\xi_n$, each particle has exactly one ancestor.
Note also that the process~$\tilde \xi_n$ was constructed in such a
way that each site can be occupied by at most $ s + 1 $
particles from $\tilde \xi_n$.
So, it follows that
process~$\tilde \xi_n$ dominates a
 Galton-Watson process
with mean offspring
being greater than or equal to
\[
 \frac{d-1-s}{d}(\bbE\eta^{(s)}+1)p
\]
(the ``worst case'' for a particle from~$\tilde\xi_n$ is when it
shares its site with another~$s$ particles from~$\tilde\xi_n$,
and all those particles have already jumped to the different
sites of level $n+1$).
Therefore, since $ p > (\bbE\eta^{(s)}+1)^{-1} $, choosing~$ d $
sufficiently large, one guarantees the
survival of the process~$\tilde \xi_n$.
This concludes the proof of Theorem~\ref{Tklim}. \qed

\medskip

Theorem~\ref{Zdlim} is a consequence of the following lemma.
\begin{lem}
\label{lemma}
Denote
\[
 \KK : = \{ x \in \bbZ^d:
      \max_{1 \le i \le d} |x^{(i)}| \le 1 \},
\]
and consider \FM{\bbZ^d}{p}{\eta}, where
$ p > {(1+\bbE\eta)}^{-1} $, {\it restricted on\/} $\KK$
(this means that if a particle attempts to jump outside $\KK$,
then it disappears).
There are constants $ d_0$, $ a > 0 $ and
$ \mu > 1 $ such that if $ d \geq d_0 $, then
with probability greater than~$ a $,
   at time $ \sqrt{d} $ there are more than $ \mu^{\sqrt d} $
active particles in~$\KK$.
\end{lem}

\medskip
\noindent
{\it Proof.} First observe that it
is enough to prove the lemma for \FM{\bbZ^d}{p}{\eta^{(s)}}
with $\eta^{(s)}$ defined by~(\ref{etas}),
where~$s$ is such that $ p > {(1+\bbE\eta^{(s)})}^{-1} $.
 Let us consider the sets
\[
 \SS := \Big\{ x \in \KK : \sum_{i=1}^d |x^{(i)}| = k \Big\},
\]
 $ k = 0, \dots, d $ and define $ \xi_k $ as the set of
active particles which are in~$\SS$ at instant $ k $.
Similarly to the proof of Theorem~\ref{Tklim}, the idea
is to show that up to time $ {\sqrt{d}} $ the process $ \xi_k $
dominates a supercritical branching process to be defined later.

Let $ x \in \SS $ and $ y \in \SSU $ be such that $ \| x-y \|=1 $,
where $ \| x\| = \sum_{i=1}^d |x^{(i)}| $.
 Notice that if site $ x $ contains an active particle at
instant~$ k $, then this particle can jump into $ y $
 at the next instant of time. Keeping
this in mind we define for $x\in\SS$
\[
 \EEE_x := \Big\{ z \in \SS :
\sum_{i=1}^d \I_{\{|x^{(i)}-z^{{(i)}}|\neq 0\}} = 2 \mbox{
and } \|x-z\| = 2 \Big\}
\]
 called the set of the ``enemies'' of $ x $. Observe that for $ x
\in \SS $ and $ z \in \EEE_x $  there exists $ y \in \SSU $
 such that $ \| x - y \| = \| z - y \| = 1 $ which in
words means that if sites $ x $ and $ z $ have active
particles at instant $ k $, then these particles can jump into the
same site next step. Moreover, for fixed~$x$ and~$z$,
the site $ y $ is the only one in $
\SSU $ with this property and there are exactly $ k + 1 $ sites in
$ \SS $ whose particles might jump into $ y $ in one step.
Notice also that $ |\EEE_x| = 2k(d-k) $.

Let
\[
 \DD_x := \{ y \in \SSU : \|x-y\| = 1 \}
\]
 be the set of ``descendants'' of~$ x \in \SS $.
It is a fact that $ |\DD_x| = 2(d-k) $. Finally, we define
for $ y \in \SSU $
\[
 \AA_y := \{ x \in  \SS : \|x-y\|=1 \},
\]
called the set of ``ancestors'' of $ y $, observing that for $ x \in \SS $
\begin{equation}
\label{disjunion}
 \EEE_x = \bigcup_{y \in \DD_x} (\AA_y \setminus \{ x \})
  \qquad \mbox{is a disjoint union,}
\end{equation}
and $|\AA_y|=k+1$ for any $ y \in \SSU $.

Note that a single site $x\in\SS$ can contain
various particles from $\xi_k$. Now (as in the proof of
Theorem~\ref{Tklim}) we define a process
$\tilde\xi_k \subset \xi_k$ in the following way.
 First, initially
the particle(s) in ${\bf 0}$ belong(s) to $\tilde\xi_0$.
If at time $k$ the set of particles $\tilde\xi_k$
(which live in $\SS$)
was constructed,
then at time~$k+1$  the set of particles $\tilde\xi_{k+1}$
(which live in $\SSU$)
is constructed in the following way.
Introduce some ordering of the particles of~$\tilde\xi_k$,
they will be allowed to jump according to that order.
Now, if the current particle survives, then
\begin{itemize}
\item if the particle jumps to some site
of~$\SSU$ and does not encounter
any particles that already belong to $\tilde\xi_{k+1}$
 there, then this particle as well as all the particles possibly
activated by it enter to $\tilde\xi_{k+1}$;
\item otherwise it is deleted.
\end{itemize}
 For $x\in\SS$ define $\XX(x)$ to be the number of particles from
$\tilde\xi_k$ in the site~$x$. Note that, by construction,
$0\leq \XX(x) \leq s+1$ for all~$x$ and~$k$.
 For $ x \in \SS $ and $ y \in \DD_x $ we denote
by $ (x \to y) $ the event
\[
 \left\{
   \begin{array}{c}
    \XX(x)\geq 1
    \mbox{ and at least one particle}\\
    \mbox{from $\tilde\xi_k$ jumps from $x$ to $y$ at time $k+1$}
   \end{array}
   \right\},
\]
and let $  \zeta^k_{x y} $
be the indicator function of the event
\[
\{\mbox{there is $z \in \EEE_x$ such that $\XX(z)\geq 1$ and $(z \to y)$}
   \}.
\]
Picking $ k \leq {\sqrt d} $, it is true that
\[
 \bbP[\zeta^k_{xy} = 1] \le
    \bbP^*[\zeta^k_{xy} = 1]
\le \frac{C_1 k}{d} \leq C_1({\sqrt d})^{-1}
\]
for some positive constant $C_1 = C_1(s)$,
where
\[
\bbP^*[~\cdot~] = \bbP[~\cdot \mid \XX(z)=s+1 \mbox{ for all }
        z\in\EEE_x].
\]
 So, given an arbitrary $\sigma > 0$, it is possible to choose~$d$
so large that $ \bbP^*[\zeta^k_{xy}=1] < \sigma $
for $ k \leq {\sqrt d} $. With this choice
for $ d $, if $ \zeta^k_x $ is the indicator function of the event
\[
\Big\{ |\DD_x \cap \{\mbox{$y\in\SSU :{}$ there exists $z\in\SS\setminus\{x\}$
  such that $(z\to y)$}\}|
> 2\sigma d \Big\},
\]
 then it follows that
\[
 \bbP[\zeta^k_x = 1 ] = \bbP\Big[ \sum_{y \in \DD_x}
\zeta^k_{xy} > 2 \sigma d \Big].
\]
 Notice that by~(\ref{disjunion}) the random variables $
\{ \zeta^k_{xy}: y \in \DD_x \} $ are independent
with respect to~$\bbP^*$. Therefore, by Lemma~\ref{Shiryaev},
we get for $ k \leq {\sqrt d} $
\begin{eqnarray}
 \bbP[ \zeta^k_x = 1] & \leq & \bbP^*\Big[ \sum_{y \in \DD_x}
\zeta^k_{xy} > 2 \sigma d \Big] \nonumber\\
& = & \bbP^*\Big[\frac{\sum_{y \in \DD_x} \zeta^k_{xy}}{2(d-k)} >
\frac{2d\sigma}{2(d-k)}\Big] \nonumber\\
& \leq &
\bbP^*\Big[\frac{\sum_{y \in \DD_x} \zeta^k_{xy}}{2(d-k)} > \sigma\Big]
\nonumber\\
 & \leq &
\exp\{-2C_2(d-k)\} \nonumber\\
 & \leq & \exp\{-C_3 d\},\label{bnm}
\end{eqnarray}
with some positive constants $ C_2, C_3$, which
depend only on~$\sigma$. Let us define the following event
\[
 B := \bigcup^{\sqrt d}_{k=1}
\bigcup_{x \in \SS} \{ \zeta^k_x = 1 \}.
\]
 Since $\eta^{(s)} \le s $ we have that $ |\tilde\xi_k| \le
(s+1)^{k+1}$. Therefore, from~(\ref{bnm}) it follows that
\[
 \bbP[B] \le {\sqrt d} \times (s+1)^{{\sqrt d}+1} \exp\{-C_3 d \},
\]
 and, as a consequence, $ \bbP[B] $ can be made arbitrarily
small for fixed $ \sigma $ and $ d $ large enough.

Suppose that the event $ B^c $ happens. In this case,
since each site can be occupied by at most $ s + 1 $ particles
 from~$\tilde\xi_k$, for each $ x
\in \SS $ there are at least $ 2(d - {\sqrt d}) - 2 \sigma d - s $
{\it available sites} (i.e., sites which do not yet contain
any particle from $\tilde\xi_{k+1}$)
in $\SSU$ into which a particle from~$\tilde\xi_k$
placed at site $ x $ could jump.
 So, it follows that
up to time~$\sqrt d$, the process $\tilde\xi_k$ dominates
a Galton-Watson branching process with mean
 offspring being greater than or equal to
\begin{equation}
\label{qwe}
 \frac{(2(d-\sqrt{d})-2 \sigma d - s)(\bbE\eta^{(s)}+1)p}{2d}.
\end{equation}
Pick $ \sigma $ small enough and $ d $ large enough to
make~(\ref{qwe}) greater than~$1$.
 The lemma follows since with positive probability a
supercritical branching process grows exponentially in time. \qed

\medskip
\noindent
{\it Proof of Theorem~\ref{Zdlim}.}
Let us first introduce some notation. Remember that
\[
 \Lambda : =
\{ x \in \bbZ^d: x^{(i)}=0 \mbox{ for } i \ge 3 \} .
\]
\noindent is a copy of $\bbZ^2$ immersed in $\bbZ^d.$
For $M\in\bbN$ denote by
\[
 \Lambda_M=\{x\in\Lambda :
\max(|x^{(1)}|,|x^{(2)}|)\leq M\} \subset \Lambda
\]
\noindent the square centered at the origin and with
sides of size $2M$, parallel to the coordinate axes.
 For $x\in \Lambda_M$ let
\[
\ell_x = \{y\in\bbZ^d : y^{(1)}=x^{(1)}, y^{(2)}=x^{(2)}\} \subset \bbZ^d
\]
\noindent be the line orthogonal to $\Lambda$
containing the site $x \in \Lambda_M$.
 By Lemma~\ref{lemma}, for $ d \geq d_0 $,
 at instant $ {\sqrt d} $ there are more than $ \mu^{\sqrt d} $
active particles, where $\mu > 1$, in $ \KK = \{ x \in \bbZ^d:
\max_{1 \le i \le d} |x^{(i)}| \le 1 \} $ with probability larger
than $ a > 0 $ for the process restricted on~$\KK$.

Fixed $ M \in \bbN, \ y \in \KK $ and $ x \in \Lambda_M$,
 after at most $ 2M+2 $
steps an active particle starting from $ y $ hits $ \ell_x $
with probability at least $({p}/{2d})^{2M+2}$.
So, for each fixed
site $ x $ of $ \Lambda_M $, the probability that
at least one of those $\mu^{\sqrt d}$ particles
 enters $\ell_x$ after at most $ 2M + 2 $
steps is greater than
\[
 1 - \Big(1 -  \Big(\frac{p}{2d}\Big)^{2M+2}
                 \Big)^{\mu^{\sqrt d}}.
\]
Consequently, defining
\begin{eqnarray*}
a' &:=& \bbP [\mbox{$\ell_x$ is hit by some
   particle starting from $\KK$,}\\
 &&~~~~~\mbox{for all $x\in\Lambda_M$}\mid \mbox{more than $\mu^{\sqrt d}$
                      particles start from $\KK$}]
\end{eqnarray*}
 one gets that, for fixed~$M$,
\[
a' \geq 1 - (2M+1)^2\Big(1-\Big(\frac{p}{2d}\Big)^{2M+2}
                 \Big)^{\mu^{\sqrt d}}
\]
and so~$a'$ can be made arbitrarily
close to~$1$ by choosing~$d$ large enough.
So we see that with probability at least $aa'\bbP[\eta\geq 1]$
the projection of the trajectories of particles from~$\KK$
will fill up the square $\Lambda_M$ (and, by choosing~$d$ large enough,
$M$ can be made as large as we want).

Note that we can repeat the above construction for any site
$x\in 3\Lambda$, and note also that if $x,y\in 3\Lambda$, $x\neq y$,
then those constructions starting from~$x$ and~$y$ are
independent (since $(\KK+x)\cap(\KK+y) = \emptyset$).
Consider now the following percolation model: For $x\in 3\Lambda$,
all the sites of the square $\Lambda_M+x$ are selected with
probability $aa'\bbP[\eta\geq 1]$. Then, as in Theorem~\ref{Zd},
one can prove that for~$M$ large enough this model percolates.
Using Lemma~\ref{H}, we obtain that the original frog model
survives with positive probability, thus concluding the
proof of Theorem~\ref{Zdlim}.
\qed

\subsection{Recurrence and transience}
\label{s_rectrans}
{\it Proof of Theorems~\ref{Tktrans}, \ref{Zdtrans},
                                      and~\ref{chi_trans}.}
The idea of the proof of all the theorems about transience in
this section is the following: {\it all\/} the particles are made
active initially; clearly, if in such model with probability~$1$
the origin is hit only a finite number of times, then a
coupling argument shows that the original frog model
is transient.

To prove Theorem~\ref{Tktrans}, we need
an upper bound for $\bbP[y\in\RR^1_x]$
which is better than~(\ref{ur*}). Note that on $\TT$, the
probability that a SRW (no death)
starting from~$x$ will eventually hit~$y$,
is exactly $(d-1)^{-\dist(x,y)}$. This shows that, on~$\TT$,
\begin{equation}
\label{ur*'}
\bbP[y\in\RR^1_x] \leq \sum_{i=\dist(x,y)}^\infty p^i(1-p)
  \frac{1}{(d-1)^{\dist(x,y)}} =
  \Big(\frac{p}{d-1}\Big)^{\dist(x,y)}.
\end{equation}
Now, using~(\ref{ur*'}) and Lemma~\ref{analiz}, one gets
that for some $C>0$
\begin{eqnarray*}
\sum_{x\neq {\bf 0}} \bbP[{\bf 0}\in \RR_x] &=&
   \sum_{i=1}^\infty \gamma_i \sum_{k=1}^\infty d(d-1)^{k-1}
      \Phi(i,k,p/(d-1)) \\
      &\leq & C\sum_{i=1}^\infty \gamma_i
       i^{\frac{\log(d-1)}{\log((d-1)/p)}} < \infty
\end{eqnarray*}
for any $p<1$, so from Borel-Cantelli one gets that almost surely
only a finite number of particles will ever enter~${\bf 0}$, thus
proving Theorem~\ref{Tktrans}.

As for Theorem~\ref{Zdtrans}, we have, recalling the proof of
Theorem~\ref{Zdext}, that when $\bbE(\log(\eta\vee 1))^d<\infty$,
\[
\sum_{x\neq {\bf 0}} \bbP[{\bf 0}\in \RR_x] =
 \sum_{x\neq {\bf 0}} \bbP[x\in \RR_{\bf 0}] =
  \bbE|\RR_{\bf 0}\setminus\{{\bf 0}\}| < \infty,
\]
and Theorem~\ref{Zdtrans} follows from Borel-Cantelli as well.

Let us turn to the proof of Theorem~\ref{chi_trans}. Denote
by~$r_k(\GG)$ the expected size of the range of the SRW on~$\GG$ until the
moment~$k$. We have
\begin{eqnarray*}
\sum_{x\neq {\bf 0}} \bbP[{\bf 0}\in \RR_x] &=&
   \sum_{x\neq {\bf 0}} \bbP[x\in \RR_{\bf 0}]\\
   &=& \bbE|\RR_{\bf 0}|\\
   &=& \sum_{k=1}^\infty \bbP[\Xi=k]r_k(\GG),
\end{eqnarray*}
and using the fact that
\[
 r_k(\GG) \simeq \left\{\begin{array}{ll}
                        \sqrt {k}, & \GG=\bbZ,\\
                        \displaystyle\frac{k}{\log k}, & \GG=\bbZ^2,\\
                        k, & \GG=\bbZ^d \mbox{ or } \TT, d\geq 3
                        \end{array}
                 \right.
\]
(see e.g.\ Hughes~\cite{Hughes},  p.\ 333, 338), one gets the result.
\qed

\medskip
\noindent
{\it Proof of Theorems~\ref{Tkrec}, \ref{Zdrec}, and~\ref{chi_rec}.}
In this section, theorems concerning the recurrence also are proved
using a common approach.
This approach can be roughly described as follows.
 We think of $\GG$ as a
disjoint union of sets $ \II_k, k=1,2,\ldots$, of
increasing sizes,
such that with large probability
(increasing with~$k$), the set $\II_k$
contains a lot of particles in the initial configuration.
Besides, given that~$\II_k$ contains many particles, also with large
probability (increasing with~$k$ as well),
the virtual paths of those particles
will cover the whole set $\II_{k+1}$ together with the origin,
 thus activating all particles placed
originally in $\II_{k+1}$. With a particular choice of
that sequence of sets, all the
events mentioned above occur simultaneously with strictly
positive probability, which implies, consequently, that the
process is recurrent (as in this case for each~$k$ there is a
particle from~$\II_k$ which visits the origin, and so the total number
of particles visiting the origin is infinite).

First, we give the proof of Theorem~\ref{Tkrec}. Fix a number~$\alpha>1$
in such a way that $\frac{\log(d-1)}{2\log(\alpha d)}>\beta$, and fix
the survival parameter~$p$ in such a way that $1/\alpha<p<1$.
Denote $\II^d_n = \{y\in\TT : \dist({\bf 0}, y)=n\}$, and define the events
\begin{eqnarray}
A_n^d &=& \{\mbox{there exists $x\in\II_n^d$ such that }
                           \eta(x)\geq (\alpha d)^{2n}\},\nonumber\\
B_n^d &=& \Big\{(\II_n^d\cup\{{\bf 0}\}) \subset
                 \bigcup_{y\in \II^d_{n-1}}\RR_y\Big\}. \label{def_B_n}
\end{eqnarray}
We have, as $|\II_n^d|>(d-1)^n$ (note also that
$\frac{d-1}{(\alpha d)^{2\beta}}>1$), that
\begin{eqnarray}
\bbP[A_n^d] &\geq & 1-(1-\bbP[\eta\geq (\alpha d)^{2n}])^{(d-1)^n}\nonumber\\
 &\geq & 1 - \Big(1-\frac{1}{(\alpha d)^{2\beta n}}\Big)^{(d-1)^n}\nonumber\\
 &\geq & 1 - C_1\exp\Big(-
   \Big(\frac{d-1}{(\alpha d)^{2\beta}}\Big)^n\Big)\label{uravn1}
\end{eqnarray}
for some $C_1>0$.
Now, using the fact that
\[
\max_{\substack{x\in\II_n^d,\\
 y\in\II_{n+1}^d\cup\{{\bf 0}\}}}\dist(x,y) = 2n+1
\]
together with~(\ref{ur**}), one gets (note that
$|\II_{n+1}^d\cup\{{\bf 0}\}|\leq d^{n+1}$ for all~$n$)
\begin{eqnarray*}
\bbP[B_{n+1}^d \mid A_n^d, B_n^d] &\geq & 1 - |\II_{n+1}^d\cup\{{\bf 0}\}|
               \Big(1-\Big(\frac{p}{d}\Big)^{2n+1}\Big)^{(\alpha d)^{2n}}\\
&\geq & 1 - C_2 d^{n+1} \exp(-(\alpha p)^n)
\end{eqnarray*}
for some $C_2>0$. The fact that $\alpha p > 1$ together with~(\ref{uravn1})
imply that with strictly positive probability there exists a
random number~$n_0$ such that the events $B_n^d, n\geq n_0$, occur.
Clearly, in this case ${\bf 0}$ is hit infinitely often and so the process
is recurrent.

Now, we start proving Theorem~\ref{Zdrec}. Fix any $0<p\leq 1$ and
choose~$\alpha>1$ in such a way that $d-\alpha\beta>0$. Let
$\II_n^d = \{x\in\bbZ^d : 2^{n-1}<\dist({\bf 0}, x)\leq 2^n\}$. Define the
events $B_n^d$ by means of~(\ref{def_B_n}) and
\[
A_n^d = \{\mbox{there exists $x\in\II_n^d$ such that }
                           \eta(x)\geq \exp(2^{\alpha n})\}.
\]
As $|\II_n^d|\geq C_1 2^{dn}$ for some $C_1>0$ and all~$n$, we have
\begin{eqnarray}
\bbP[A_n^d] &\geq & 1-(1-
   \bbP[\eta\geq \exp(2^{\alpha n})])^{C_1 2^{dn}}\nonumber\\
 &\geq & 1 - \Big(1-\frac{1}{2^{\alpha\beta n}}\Big)^{C_1 2^{dn}}\nonumber\\
 &\geq& 1 - C_2 \exp(-(C_1 2^{d-\alpha\beta})^n).\label{uravn2}
\end{eqnarray}
It is a fact that in this case
\[
\max_{\substack{x\in\II_n^d,\\ y\in\II_{n+1}^d\cup\{{\bf 0}\}}}\dist(x,y)
   \leq 2^{n+2},
\]
and that $|\II_{n+1}^d\cup\{{\bf 0}\}|\leq C_3 2^{d(n+1)}$,
so using~(\ref{ur**}) we get, keeping in mind that $\alpha>1$,
\begin{eqnarray}
\bbP[B_{n+1}^d \mid A_n^d, B_n^d] &\geq & 1 - |\II_{n+1}^d\cup\{{\bf 0}\}|
               \Big(1-\Big(\frac{p}{2d}\Big)^{2^{n+2}}
                  \Big)^{\exp(2^{\alpha n})}\nonumber\\
&\geq & 1 - C_4 2^{d(n+1)}\exp\Big(-\exp\Big(2^{\alpha n} -
               2^{n+2}\log\frac{2d}{p}\Big)\Big).
                 \phantom{****} \label{uravn3}
\end{eqnarray}
As before, (\ref{uravn2})--(\ref{uravn3}) imply that with positive
probability infinite number of events $B_n^d$ occur, so the process
is recurrent.

Let us turn to the proof Theorem~\ref{chi_rec}.
The sets $\II_n^d$ are now defined by
$\II_n^d = \{x\in\bbZ^d : 2^{n-1}<\|x\|_2\leq 2^n\}$, and
the sequence of events~$B_n^d$ is still defined by~(\ref{def_B_n}).
Recall that $\Xi^x$ is the lifetime of the particle originating
from~$x$. Now, the site $x\in\II_n^d$ is called {\it good}, if
$\Xi^x\geq 2^{2(n+2)}$ (intuitively, the site~$x\in\II_n^d$
is good if the corresponding particle lives long enough
to be able to get to any fixed site
of $\II_{n+1}^d\cup\{{\bf 0}\}$). Define the events
\[
 A_n^d = \{\mbox{the number of good sites in }\II_n^d
            \geq \phi_{n,d}|\II_n^d|\},
\]
where
\[
\phi_{n,d} = \left\{
             \begin{array}{ll}
              \beta_1 2^{-(n+3)} \log\log 2^{n+2}, & d=1, \\
              \beta_2 2^{-(2n+5)} (\log 2^{n+2})^2, & d=2,
                           \phantom{\displaystyle\sum^a} \\
              \beta_d 2^{-(2n+5)} \log 2^{n+2}, & d\geq 3
                           \phantom{\displaystyle\sum^a}. \\
             \end{array}
             \right.
\]
As, by the hypothesis, $\bbP[x \mbox{ is good}]\geq 2\phi_{n,d}$
for any~$x\in\II_n^d$, by Lemma~\ref{Shiryaev} we get
(observe that $|\II_n^d|\simeq 2^{dn}$)
\[
\bbP[A_n^d] \geq \left\{
             \begin{array}{ll}
              1-n^{-C_1\beta_1}, & d=1,\\
              1-\exp(-C_2\beta_2n^2), & d=2,
                \phantom{\displaystyle\sum^a}\\
              1-\exp(-C_3\beta_d n 2^{d(n-2)}), & d\geq 3,
                 \phantom{\displaystyle\sum^a}
             \end{array}
             \right.
\]
so $\sum_{n=1}^\infty (1-\bbP[A_n^d]) < \infty$ for any $\beta_d$,
$d\geq 2$, and for $\beta_1>1/C_1$, $d=1$.
Using the inequality
\[
\max_{\substack{x\in\II_n^d,\\ y\in\II_{n+1}^d\cup
\{{\bf 0}\}}}\|x-y\|_2
   \leq 2^{n+2}
\]
together with Lemma~\ref{pofvisit}, one gets
\[
\bbP[B_{n+1}^d \mid A_n^d, B_n^d] \geq \left\{
             \begin{array}{ll}
              1-3(1-w_1)^{\phi_{n,1}|\II_n^1|}, &
               d=1,\\
              1-|\II_{n+1}^2\cup\{{\bf 0}\}|\Big(1-
               \displaystyle\frac{w_2}{\log 2^{n+2}}\Big)^{\phi_{n,2}
                |\II_n^2|}, &
                d=2,\\
              1-|\II_{n+1}^d\cup\{{\bf 0}\}|\Big(1-
               \displaystyle\frac{w_d}{2^{(d-2)(n+2)}}\Big)^{\phi_{n,d}
                |\II_n^d|}, &
                d\geq 3
             \end{array}
             \right.
\]
(the factor $|\II_{n+1}^d\cup\{{\bf 0}\}|$ was substituted by~3
in dimension~1, because in this case, to guarantee that all
the sites of the set $\II_{n+1}^1\cup\{{\bf 0}\}$ are hit, it is
sufficient to visit the sites~${\bf 0}$ and~$\pm 2^{n+1}$).
Then, elementary computations show that
\[
\bbP[B_{n+1}^d \mid A_n^d, B_n^d] \geq \left\{
             \begin{array}{ll}
              1-3n^{-C_4\beta_1}, & d=1,\\
              1-C_5 2^{-(C_6w_2\beta_2-2)n}, & d=2,
                   \phantom{\displaystyle\sum^a}\\
              1-C_7 2^{-(C_8w_d\beta_d-d)n}, & d\geq 3.
                \phantom{\displaystyle\sum^a}
             \end{array}
             \right.
\]
By choosing $\beta_1 > \max\{1/C_1,1/C_4\}$, $\beta_2>2/C_6w_2$,
$\beta_d>d/C_8w_d$, $d\geq 3$, once again one gets that with positive
probability infinite number of events~$B_n^d$ occur,
and so the frog model is recurrent. \qed

\section*{Acknowledgements}
This work was done when O.S.M.~Alves was visiting the Statistics
Department of the Institute of Mathematics and Statistics of the
University of S\~ao Paulo. He is thankful to the probability group
of IME-USP for hospitality.
Special thanks are due to the
anonymous referee of the first version of this paper, who
suggested the method of proving Theorems~\ref{Tkext},
\ref{Zdext}, and~\ref{Tksurv}.


\begin{thebibliography}{9}

\bibitem{sofa} {\sc O.S.M.~Alves, F.P.~Machado, S.Yu.~Popov} (2000)
{\it The shape theorem for the frog model.} (submitted).
This article is available from the xxx mathematics archive as math.PR/0102182
\\
(http://front.math.ucdavis.edu/math.PR/0102182).

%\bibitem{Bill} {\sc P.~Billingsley} (1995)
%{\it Probability and Measure} (3rd.\ ed.), John Wiley and Sons,
%Inc., New York.

%\bibitem{H} {\sc O.~H\"aggstr\"om} (2000)
% Markov random fields and percolation on general graphs.
% {\it Adv. Appl. Probab.} {\bf 32} (1), 39--66.

\bibitem{Ravi} {\sc O.S.M.~Alves, F.P.~Machado, S.Yu.~Popov,
         K.~Ravishankar} (2001)
{\it The shape theorem for the frog model with
random initial configuration.} (in preparation).

\bibitem{HMP} {\sc F.~den~Hollander, M.V. Menshikov, S.Yu. Popov} (1999)
   A note on transience versus recurrence for a branching random walk
    in random environment.
         {\it J. Statist. Phys.} {\bf 95} (3/4), 587--614.

\bibitem{Hughes} {\sc B.D.~Hughes} (1995)
{\it Random Walks and Random Environments, vol.~1.}
Clarendon press, Oxford.

\bibitem{Lawler} {\sc G.F.~Lawler} (1991)
{\it Intersections of Random Walks.} Birkh\"auser Bos\-ton.

\bibitem{M} {\sc M.V.~Menshikov} (1985)
Estimates for percolation thresholds for lattices in $\bbR^d$.
{\it Soviet Math. Dokl.} {\bf 32} (2), 368--370.

\bibitem{Micha} {\sc M.V.~Menshikov} (1987)
Qualitative estimates and rigorous inequalities for critical
points of a graph and its subgraphs. {\it Probab. Theory Appl.}
{\bf 32}, 544--547.

\bibitem{Serguei} {\sc S.Yu.~Popov} (2001)
Frogs in random environment.
{\it J. Statist. Phys.} {\bf 102} (1/2), 191--201.

\bibitem{Shir} {\sc A. Shiryaev} (1989)
{\it Probability} (2nd. ed.). Springer, New York.

\bibitem{TW} {\sc A.~Telcs, N.C.~Wormald} (1999)
Branching and tree indexed random walks on fractals.
{\it J. Appl. Probab.} {\bf 36} (4), 999--1011.

\end{thebibliography}
\end{document}